# VARIANCE-OPTIMAL HEDGING FOR PROCESSES WITH STATIONARY INDEPENDENT INCREMENTS[1]

By Friedrich Hubalek, Jan Kallsen and Leszek Krawczyk

*University of Aarhus, Munich University of Technology and Vienna University of Technology*

We determine the variance-optimal hedge when the logarithm of the underlying price follows a process with stationary independent increments in discrete or continuous time. Although the general solution to this problem is known as backward recursion or backward stochastic differential equation, we show that for this class of processes the optimal endowment and strategy can be expressed more explicitly. The corresponding formulas involve the moment, respectively, cumulant generating function of the underlying process and a Laplace- or Fourier-type representation of the contingent claim. An example illustrates that our formulas are fast and easy to evaluate numerically.

**1. Introduction.** A basic problem in mathematical finance is how an option writer can hedge her risk by trading only in the underlying asset. This question is well understood in frictionless complete markets. It suffices to buy the replicating portfolio in order to completely offset the risk. This elegant approach works well in the standard Black–Scholes or Cox–Ross–Rubinstein setup, but not much beyond.

On the other hand, it has often been reported that real market data exhibit heavy tails and volatility clustering. Two common ways to account for such phenomena are some sort of stochastic volatility or jump processes, or a combination of both. In this paper we adopt the second approach and assume that the logarithmic stock price follows a general process with stationary, independent increments, either in discrete or continuous time. Processes of

Received December 2004; revised October 2005.

[1]Supported in part by Austrian Science Foundation (FWF) Grant SFB 10 ("Adaptive Information Systems and Modelling in Economics and Management Science") and Project P11544.

this type play by now an important role in the modeling of financial data (cf. [3, 13, 14, 36]).

Since replicating portfolios typically do not exist in such incomplete markets, one has to choose alternative criteria for reasonable hedging strategies. If you want to be as safe as in the complete case, you should invest in a *superhedging* strategy (cf., e.g., [16]). In this case you may "suffer" profits but no losses at maturity of the derivative, which is very agreeable. On the other hand, even for simple European call options, only trivial superhedging strategies exist in a number of reasonable market models ("buy the stock"; cf. [8, 12, 22]).

Alternatively, you may maximize some expected utility among all portfolios that differ only in the underlying asset and have a fixed position in the contingent claim. Variations of this approach have been investigated in [8, 9, 19, 30, 31].

In this paper we follow a third popular suggestion, namely, to minimize some form of quadratic risk (cf. [11, 21, 50, 52] for an overview). This can be interpreted as a special case of the second approach if we allow for quadratic utility functions.

Quadratic hedging comes about in two different flavors: *local risk minimization* as in [20] and [49], and *global risk minimization* (i.e., *variance-optimal hedging*, *mean-variance hedging*) as in [11] and [50]. Roughly speaking, one may say that locally optimal strategies are relatively easy to compute but hard to interpret economically, whereas the opposite is true for the globally optimal hedge. This paper focuses on the second problem, but as a by-product, we also obtain the locally optimal Föllmer–Schweizer hedge. In discounted terms, the global problem can be stated as follows: If $H$ denotes the payoff of the option and $S$ denotes the underlying price process, try to minimize the squared $L^2$-distance

$$E((c + G_T(\theta) - H)^2) \tag{1}$$

over all initial endowments $c \in \mathbb{R}$ and all in some sense admissible trading strategies $\theta$. Here, $G_T(\theta) = \int_0^T \theta_t \, dS_t$ [resp. $G_T(\theta) = \sum_{n=1}^T \theta_n \Delta S_n$ in discrete time] denotes the cumulative gains from trade up to time $T$. The idea is obviously to approximate the claim as closely as possible in an $L^2$ sense. Even though one may argue that one should not punish gains, the clarity and simplicity of this criterion are certainly appealing. Since it is harder to explain, we do not discuss local risk minimization here, but refer instead to [52].

By way of duality, quadratic optimization problems are related to (generally signed) martingale measures, namely the *Föllmer–Schweizer* or *minimal martingale measure* for local optimization, and the *variance-optimal martingale measure* for global optimization. A similar duality has been established and exploited in many recent papers on related problems of utility maximization or portfolio optimization (cf. [7, 8, 17, 18, 23, 24, 25, 26, 33, 34, 40, 47]). Roughly speaking, the minimal martingale measure is the martingale measure whose density can be written as $1 + \int_0^T \theta_t \, dM_t$ for some $\theta$, where $M$ denotes the martingale part in the Doob–Meyer decomposition of $S$. The integrand $\theta$ can be determined relatively easily in terms of the local behavior of $S$, which may be given by a stochastic differential equation or by one-step transition probabilities in discrete time. By contrast, the variance-optimal martingale measure is characterized by a density of the form $c + \int_0^T \theta'_t \, dS_t$ for some $c \in \mathbb{R}$ and some (generally different) integrand $\theta'$. Here, it is usually much harder to determine $\theta'$. This holds with one notable exception, namely if the so-called *mean-variance trade-off* process is deterministic, in which case both measures coincide. More specifically, the integrands $\theta$ and $\theta'$ above tally because the difference $\int_0^T \theta_t \, dS_t - \int_0^T \theta_t \, dM_t$ is a constant and can be moved to $c$. In this case of deterministic mean-variance trade-off, globally risk-minimizing hedging strategies can be computed from locally risk-minimizing ones. The setup in this paper is among the few models of practical importance where the condition of deterministic mean-variance trade-off naturally holds.

The process formed by conditional expectation of the option's payoff under the minimal, respectively, variance-optimal martingale measure is sometimes interpreted as a derivative price process. In jump-type models one has to be careful at this point because these measures are generally signed and may lead to arbitrage. In particular, the variance-optimal initial capital for a positive payoff may be negative; see [2], Section 4 and the numerical example at the end of Section 5.2 below.

Even in the case of deterministic mean-variance trade-off, the actual computation of variance-optimal hedging strategies involves the joint predictable covariation of the option's "price process" and the underlying stock. For general claims, how to obtain this covariation may not seem evident. It can be computed quite easily if the payoff is of exponential type $e^{zX_T}$, where $X := \log(S/S_0)$ denotes the process with stationary, independent increments driving the stock price $S$. The reason is that the "price process" for such exponential payoffs under the variance-optimal martingale measure is again the exponential of a process with stationary independent increments, which leads to handy formulas for the corresponding hedge. Since the optimality criterion in (1) is based on an $L^2$-distance, the resulting hedging strategy is linear in the option. This suggests writing an arbitrary claim as a linear combination of exponential payoffs. Put differently, we work with the inverse Laplace (or Fourier) transform of the option. This will be done in Section 2 for discrete-time and in Section 3 for continuous-time processes, respectively. One could go even one step further and generalize the results to arbitrary processes with independent increments because they still share the important property of deterministic mean-variance trade-off. However, we chose not to do so in order not to drown the arguments in technicalities and because this more general class plays a minor role in applications.

Since the first version of this paper circulated, the idea of using Fourier or Laplace transforms with Lévy processes has been applied independently in the framework of option pricing by [4] as well as [42] and, very recently, in the context of quadratic hedging by [53]. Explicit integral representations of a number of concrete payoffs are to be found in Section 4.

Section 5 illustrates the application of the results. We compare the variance-optimal hedge of a European call in a pure-jump Lévy process model to the Black–Scholes hedge as a benchmark. Since the results in the subsequent sections rely heavily on bilateral Laplace transforms, the Appendix contains a short summary of important results in this context.

To keep the presentation and notation simple, we confine ourselves to one single underlying asset. Extensions to the multivariate case and to path-dependent claims will be provided elsewhere. For unexplained notation, we refer the reader to standard textbooks on stochastic calculus such as, for example, [29] or [41].

**2. Discrete time.** Let $(\Omega, \mathcal{F}, (\mathcal{F}_n)_{n \in \{0,1,\ldots,N\}}, P)$ denote a filtered probability space and let $X = (X_n)_{n=0,1,\ldots,N}$ denote a real-valued process with stationary, independent increments in the sense that:

1. $X$ is adapted to the filtration $(\mathcal{F}_n)_{n \in \{0,1,\ldots,N\}}$;
2. $X_0 = 0$;
3. $\Delta X_n := X_n - X_{n-1}$ has the same distribution for $n = 1, \ldots, N$;
4. $\Delta X_n$ is independent of $\mathcal{F}_{n-1}$ for $n = 1, \ldots, N$.

We consider a nondividend-paying stock whose discounted price process $S$ is of the form

$$S_n = S_0 \exp(X_n)$$

with some constant $S_0 > 0$. We assume that $E(S_1^2) = S_0^2 E(e^{2X_1}) < \infty$, which implies that the moment generating function $m : z \mapsto E(e^{zX_1})$ is defined at least for $z \in \mathbb{C}$ with $0 \leq \Re(z) \leq 2$. Moreover, we exclude the degenerate case that $S$ is deterministic. Put differently, $\mathrm{Var}(e^{X_1}) = m(2) - m(1)^2$ does not vanish, which can be viewed as a no arbitrage condition.

Our goal is to determine the variance-optimal hedge for a European-style contingent claim on the stock that expires at $N$ with discounted payoff $H$. Mathematically, $H$ denotes a square-integrable, $\mathcal{F}_N$-measurable random variable of the form $H = f(S_N)$ for some function $f : (0, \infty) \to \mathbb{R}$. More specifically, we assume that $f$ is of the form

(2) $$f(s) = \int s^z \Pi(dz)$$

for some finite complex measure $\Pi$ on a strip $\{z \in \mathbb{C} : R' \leq \Re(z) \leq R\}$, where $R', R \in \mathbb{R}$ are chosen such that $E(e^{2R'X_1}) < \infty$ and $E(e^{2RX_1}) < \infty$. The integral representation of European call and put options, and several other payoffs can be found in Section 4.

REMARK 2.1. Loosely speaking, the option's payoff at maturity is written as a linear combination of powers of the $S$ or exponentials of $X$. Put differently, its payoff function is a kind of inverse Mellin or Laplace transform of the measure $\Pi$. To be more specific, let us consider the case $R' = R$, that is, $\Pi$ is concentrated on the straight line $R + i\mathbb{R}$. Denote by $\ell$ the inverse Laplace transform of $\Pi$ in the sense that $\ell(x) = \int_{R-i\infty}^{R+i\infty} e^{zx} \Pi(dz)$ for $x \in \mathbb{R}$. Then

$$H = f(S_N) = f(\exp(X_N + \log(S_0))) = \ell(X_N + \log(S_0)).$$

Up to a factor $e^{Rx}$, the function $\ell$ is just the Fourier transform of a complex measure on the real line [namely the measure $\nu$ with $\nu(B) = \Pi(R + iB)$ for Borel sets $B \subset \mathbb{R}$]. The reason to allow also $R \neq 0$ is simply that $\ell$ cannot be written as the Fourier transform of a finite measure for important claims as, for example, European calls. Examples with $R' < R$ are given in Remark 4.1 part 2.

The variance-optimal hedge minimizes the $L^2$-distance between the option's payoff and the terminal value of the hedging portfolio. To be more specific, define the set $\Theta$ of *admissible* strategies as the set of all predictable processes $\theta$ such that the *cumulative gains* $G_n(\theta) := \sum_{k=1}^n \theta_k \Delta S_k$ are square-integrable for $n = 1, \ldots, N$. We call $\phi \in \Theta$ a *variance-optimal hedging strategy* and $V_0 \in \mathbb{R}$ a *variance-optimal initial capital* if $c = V_0$ and $\theta = \phi$ minimize the expected squared hedging error

(3) $$E((c + G_N(\theta) - H)^2)$$

over all initial endowments $c \in \mathbb{R}$ and all admissible strategies $\theta \in \Theta$. Let us emphasize that the variance-optimal initial capital is, in general, not an arbitrage-free price and can be negative for a positive payoff; see Section 5.2 for a concrete example.

In our framework the variance-optimal hedge and its corresponding hedging error can be determined quite explicitly:

THEOREM 2.1. *The variance-optimal initial capital $V_0$ and the variance-optimal hedging strategy $\phi$ are given by*

(4) $$V_0 = H_0$$

*and the recursive expression*

(5) $$\phi_n = \xi_n + \frac{\lambda}{S_{n-1}}(H_{n-1} - V_0 - G_{n-1}(\phi)),$$

*where the processes $H_n$ and $\xi_n$ and the constant $\lambda$ are defined by*

(6)
$$g(z) := \frac{m(z+1) - m(1)m(z)}{m(2) - m(1)^2},$$
$$h(z) := m(z) - (m(1) - 1)g(z),$$
$$\lambda := \frac{m(1) - 1}{m(2) - 2m(1) + 1},$$
$$H_n := \int S_n^z h(z)^{N-n} \Pi(dz),$$
$$\xi_n := \int S_{n-1}^{z-1} g(z) h(z)^{N-n} \Pi(dz).$$

*The optimal hedge $V_0, \phi$ is unique up to a null set.*

REMARK 2.2. One may also consider a similar problem where the initial endowment is fixed and the mean squared difference in (3) is minimized only over the strategies $\theta \in \Theta$. This *risk-minimizing hedging strategy* for given initial capital $c$ is determined as in Theorem 2.1 with $V_0 = c$ instead of (4).

THEOREM 2.2. *The variance of the hedging error $E((V_0 + G_N(\phi) - H)^2)$ in Theorem 2.1 equals*
$$J_0 := \int \int J_0(y, z) \Pi(dy) \Pi(dz),$$

*where*
$$a(y, z) := h(y)h(z) \frac{m(2) - m(1)^2}{m(2) - 2m(1) + 1},$$
$$b(y, z) := m(y + z)$$
$$- (m(2)m(y)m(z) - m(1)m(y+1)m(z)$$
$$- m(1)m(y)m(z+1) + m(y+1)m(z+1))(m(2) - m(1)^2)^{-1},$$
$$J_0(y, z) := \begin{cases} S_0^{y+z} b(y, z) \dfrac{a(y,z)^N - m(y+z)^N}{a(y,z) - m(y+z)}, & \text{if } a(y,z) \neq m(y+z), \\ S_0^{y+z} b(y, z) N m(y+z)^{N-1}, & \text{if } a(y,z) = m(y+z). \end{cases}$$

The remainder of this subsection is devoted to the proofs of Theorems 2.1 and 2.2. As noted by [51], the variance-optimal hedge can be obtained from the option's Föllmer–Schweizer decomposition if the so-called mean-variance trade-off process of the option is deterministic. The latter is defined as
$$K_n := \sum_{k=1}^{n} \frac{(E(\Delta S_k | \mathcal{F}_{k-1}))^2}{\operatorname{Var}(\Delta S_k | \mathcal{F}_{k-1})} = \frac{(m(1)-1)^2}{m(2) - m(1)^2} n.$$

The Föllmer–Schweizer decomposition plays a key role in the determination of locally risk-minimizing strategies in the sense of [20] and [49], and it is defined as follows.

DEFINITION 2.1. Denote by $S = S_0 + M + A$ the Doob decomposition of $S$ into a martingale $M$ and a predictable process $A$. The sum $H = H_0 + \sum_{n=1}^{N} \xi_n \Delta S_n + L_N$ is called a *Föllmer–Schweizer decomposition* of $H \in L^2(P)$ if $H_0$ is $\mathcal{F}_0$-measurable, $\xi \in \Theta$ and $L$ is a square-integrable martingale with $L_0 = 0$ that is orthogonal to $M$ (in the sense that $LM$ is a martingale). We will use this notion as well if $H, H_0, \xi$ and $L$ are complex-valued, in which case we require $\Re(\xi) \in \Theta$ and $\Im(\xi) \in \Theta$.

In discrete time any square-integrable random variable admits such a decomposition, which can be found by a backward recursion (cf. [51], Proposition 2.6). However, since this method does not yield a closed-form solution in our framework, we do not use these results. Instead we proceed in two steps. First, we determine the Föllmer–Schweizer decomposition for options whose payoff is of power type. Second, we consider claims which are linear combinations of such options in the sense of (2). Here, we rely on the linearity of the Föllmer–Schweizer decomposition in the claim.

LEMMA 2.1. *Let $z \in \mathbb{C}$ with $S_1^z \in L^2(P)$. Then $H(z) = S_N^z$ admits a Föllmer–Schweizer decomposition $H(z) = H(z)_0 + \sum_{n=1}^{N} \xi(z)_n \Delta S_n + L(z)_N$, where*

$$
\begin{aligned}
H(z)_n &= h(z)^{N-n} S_n^z, \\
\xi(z)_n &= g(z) h(z)^{N-n} S_{n-1}^{z-1}, \\
L(z)_n &= H(z)_n - H(z)_0 - \sum_{k=1}^{n} \xi(z)_k \Delta S_k
\end{aligned}
\tag{7}
$$

*and $g(z), h(z)$ are defined in Theorem 2.1.*

PROOF. The statement could be derived from Proposition 2.6 and Lemma 2.7 of [51], but it is easier to prove it directly. Since $S_1^z$ is square-integrable, all the involved expressions are well defined. From (7) it follows that

$$\Delta L(z)_n = S_{n-1}^z h(z)^{N-n}(e^{z \Delta X_n} - h(z) - g(z)(e^{\Delta X_n} - 1)). \tag{8}$$

Since

$$
\begin{aligned}
E(e^{z \Delta X_n} &- h(z) - g(z)(e^{\Delta X_n} - 1)) \\
&= m(z) - h(z) - g(z)(m(1) - 1) \\
&= 0,
\end{aligned}
\tag{9}
$$

this implies that $E(\Delta L(z)_n|\mathcal{F}_{n-1}) = 0$ and hence $L(z)$ is a martingale.

The Doob decomposition $S = S_0 + M + A$ of $S$ satisfies

(10) $$\Delta A_n = E(\Delta S_n|\mathcal{F}_{n-1}) = S_{n-1}(m(1) - 1)$$

and hence $\Delta M_n = S_{n-1}(e^{\Delta X_n} - m(1))$. In view of (8), we obtain

$$\Delta M_n \Delta L(z)_n = S_{n-1}^{z+1} h(z)^{N-n} (e^{\Delta X_n} - m(1))(e^{z \Delta X_n} - h(z) - g(z)(e^{\Delta X_n} - 1)).$$

From

$$\begin{aligned} E(e^{\Delta X_n}(e^{z\Delta X_n} - h(z) - g(z)(e^{\Delta X_n} - 1))) \\ = m(z+1) - h(z)m(1) - g(z)m(2) + g(z)m(1) \\ = 0 \end{aligned}$$

and (9) it follows that $E(\Delta M_n \Delta L(z)_n|\mathcal{F}_{n-1}) = 0$ and hence $ML(z)$ is a martingale as well. □

PROPOSITION 2.1. *Any contingent claim $H = f(S_N)$ as in the beginning of this subsection admits a Föllmer–Schweizer decomposition $H = H_0 + \sum_{n=1}^{N} \xi_n \Delta S_n + L_N$. Using the notation of the previous lemma, it is given by*

$$H_n = \int H(z)_n \Pi(dz),$$

$$\xi_n = \int \xi(z)_n \Pi(dz),$$

$$L_n = \int L(z)_n \Pi(dz) = H_n - H_0 - \sum_{k=1}^{n} \xi_k \Delta S_k.$$

*Moreover, the processes $H_n$, $\xi_n$ and $L_n$ are real-valued.*

PROOF. First, note that $\int E(|\Delta L(z)_n|^2)|\Pi|(dz) < \infty$, where $|\Pi|$ denotes the total variation measure of $\Pi$ in the sense of [44], Section 6.1. From Fubini's theorem we conclude that

$$\begin{aligned} E(\Delta L_n \mathbb{1}_B) &= E\left(\int \Delta L(z)_n \Pi(dz)\mathbb{1}_B\right) \\ &= \int E(\Delta L(z)_n \mathbb{1}_B) \Pi(dz) = 0 \end{aligned}$$

for any $B \in \mathcal{F}_{n-1}$. Hence $L$ is a martingale. Similarly, it is shown that $ML$ is a martingale as well. The assertion concerning the decomposition follows from Lemma 2.1.

Since $H$ and $S_n$ are real-valued, we have

$$0 = (H_0 - \overline{H}_0) + \sum_{n=1}^{N}(\xi_n - \overline{\xi}_n)\Delta S_n + (L_N - \overline{L}_N),$$

which implies that $0 = \Im(H_0) + \sum_{n=1}^{N} \Im(\xi_n)\Delta S_n + \Im(L_N)$. Since the Föllmer–Schweizer decomposition of 0 is unique (cf. [48], Remark 4.11), we have that $H_0$, $\xi_n$ and $L_n$ are real-valued for $n = 1, \ldots, N$. □

Finally, we apply the preceding results to determine the variance-optimal hedge.

PROOF OF THEOREM 2.1. As observed by [48], Proposition 5.5, the process $S$ has deterministic mean-variance trade-off. From Proposition 2.1 and [51], Theorem 4.4, it follows that the variance-optimal hedging strategy $\phi$ satisfies

$$\phi_n = \xi_n + \lambda_n(H_{n-1} - H_0 - G_{n-1}(\phi))$$

with

$$\lambda_n := \frac{\Delta A_n}{E(\Delta S_n^2 | \mathcal{F}_{n-1})} = \frac{\lambda}{S_{n-1}}$$

[cf. (10)]. Moreover, the variance-optimal initial capital equals $V_0$.

For the uniqueness statement, suppose that $\widetilde{V}_0 \in \mathbb{R}$ and $\widetilde{\phi} \in \Theta$ lead to a variance-optimal hedge as well. Define $\hat{V}_0 := \frac{1}{2}(V_0 + \widetilde{V}_0)$ and $\hat{\phi} := \frac{1}{2}(\phi + \widetilde{\phi}) \in \Theta$. It is easy to verify that we would have

$$E((\hat{V}_0 + G_N(\hat{\phi}) - H)^2) < \tfrac{1}{2}E((V_0 + G_N(\phi) - H)^2) \\ + \tfrac{1}{2}E((\widetilde{V}_0 + G_N(\widetilde{\phi}) - H)^2),$$

contradicting the optimality of $(V_0, \phi)$ and $(\widetilde{V}_0, \widetilde{\phi})$, if $V_0 + G_N(\phi)$ and $\widetilde{V}_0 + G_N(\widetilde{\phi})$ did not coincide almost surely. Hence

$$V_0 - \widetilde{V}_0 + G_N(\phi - \widetilde{\phi}) = 0.$$

In particular, $G_N(\phi - \widetilde{\phi})$ is $\mathcal{F}_{N-1}$-measurable. We obtain

$$\begin{aligned}0 &= \operatorname{Var}(G_N(\phi - \widetilde{\phi}) | \mathcal{F}_{N-1}) \\ &= \operatorname{Var}((\phi - \widetilde{\phi})_N \Delta S_N | \mathcal{F}_{N-1}) \\ &= ((\phi - \widetilde{\phi})_N S_{N-1})^2 (m(2) - m(1)^2),\end{aligned}$$

which implies that $(\phi - \widetilde{\phi})_N = 0$ almost surely. By induction, we conclude that $(\phi - \widetilde{\phi})_n = 0$ for $n = N - 1, \ldots, 1$ and hence also $V_0 = \widetilde{V}_0$.

The remark following Theorem 2.1 follows from [51], Proposition 4.3. □

PROOF OF THEOREM 2.2. According to [51], Theorem 4.4, the variance of the hedging error equals

(11) $$\sum_{n=1}^{N} E((\Delta L_n)^2) \prod_{k=n+1}^{N} (1 - \lambda_k \Delta A_k)$$

with $\lambda_k = \lambda/S_{k-1}$ and $\Delta A_k$ as in (10). Since $\Delta L_n = \int \Delta L(z)_n \Pi(dz)$, we have that

$$(\Delta L_n)^2 = \int \int \Delta L(y)_n \Delta L(z)_n \Pi(dy) \Pi(dz)$$

and hence

(12) $$E((\Delta L_n)^2) = \int \int E(\Delta L(y)_n \Delta L(z)_n) \Pi(dy) \Pi(dz)$$

by Fubini's theorem. Equation (8) implies

$$\Delta L(y)_n \Delta L(z)_n = S_{n-1}^{y+z} h(y)^{N-n} h(z)^{N-n} (e^{y \Delta X_n} - h(y) - g(y)(e^{\Delta X_n} - 1))$$
$$\times (e^{z \Delta X_n} - h(z) - g(z)(e^{\Delta X_n} - 1)).$$

Since $E(S_{n-1}^{y+z}) = S_0^{y+z} m(y+z)^{n-1}$ and so forth, we have

$$E(\Delta L(y)_n \Delta L(z)_n) = S_0^{y+z} (h(y)h(z))^{N-n} m(y+z)^{n-1} b(y,z)$$

with

$$b(y,z) = m(y+z) - m(y)(h(z) - g(z)) - m(y+1)g(z)$$
$$- (h(y) - g(y))(m(z) - h(z) + g(z) - g(z)m(1))$$
$$- g(y)(m(z+1) - m(1)(h(z) - g(z)) - g(z)m(2)).$$

This expression coincides actually with $b(y,z)$ in the statement of the theorem. Consequently, we have shown

$$\sum_{n=1}^N E(\Delta L(y)_n \Delta L(z)_n) \prod_{k=n+1}^N (1 - \lambda_k \Delta A_k)$$

$$= S_0^{y+z} b(y,z) a(y,z)^{N-1} \sum_{n=1}^N \left( \frac{m(y+z)}{a(y,z)} \right)^{n-1}$$

$$= S_0^{y+z} b(y,z) \frac{a(y,z)^N - m(y+z)^N}{a(y,z) - m(y+z)}$$

unless the denominator vanishes in the last equation. In view of (11) and (12), we are done. □

Let us briefly discuss the structure of the variance-optimal hedge. The process $\xi$ in the Föllmer–Schweizer decomposition coincides with the locally risk-minimizing strategy. The process $H_n = H_0 + \sum_{k=1}^n \xi_k \Delta S_k + L_n$ that appears on the right-hand side of the Föllmer–Schweizer decomposition may be interpreted as a "price process" of the option. However, since this process may generate arbitrage, one should be careful with this interpretation. Note, however, that the difference between the locally and globally optimal hedging strategy in (5) is proportionate to the difference between this "option price" $H_{n-1}$ and the investor's current wealth.

**3. Continuous time.** We turn now to the continuous-time counterpart of the previous section. Similarly as before, $(\Omega, \mathcal{F}, (\mathcal{F}_t)_{t \in [0,T]}, P)$ denotes a filtered probability space and $X = (X_t)_{t \in [0,T]}$ denotes a real-valued process with stationary, independent increments (*PIIS, Lévy process*) in the sense that:

1. $X$ is adapted to the filtration $(\mathcal{F}_t)_{t \in [0,T]}$ and has cadlag paths;
2. $X_0 = 0$;
3. the distribution of $X_t - X_u$ depends only on $t - u$ for $0 \le u \le t \le T$;
4. $X_t - X_u$ is independent of $\mathcal{F}_u$ for $0 \le u \le t \le T$.

As in the discrete-time case, the distribution of the whole process $X$ is determined by the law of $X_1$. The latter is an infinitely divisible distribution which can be expressed in terms of its Lévy–Khinchine representation. Alternatively, one may characterize it by its *cumulant generating function*, that is, by the continuous mapping $\kappa : D \to \mathbb{C}$ with $E(e^{zX_t}) = e^{t\kappa(z)}$ for $z \in D := \{z \in \mathbb{C} : E(e^{\Re(z)X_1}) < \infty\}$ and $t \in \mathbb{R}_+$. For details on Lévy processes and unexplained notation, we refer the reader to [29, 41, 46].

The discounted price process $S$ of the nondividend-paying stock under consideration is supposed to be of the form

$$S_t = S_0 \exp(X_t)$$

with some constant $S_0 > 0$. Again, we assume that $E(S_1^2) = S_0^2 E(e^{2X_1}) < \infty$, which means that $z \in D$ for any complex number $z$ with $0 \le \Re(z) \le 2$. Moreover, we exclude the degenerate case that $S$ is deterministic, that is, we have $\kappa(2) - 2\kappa(1) \ne 0$. This can be viewed as a no arbitrage condition.

As in Section 2 we consider an option with discounted payoff $H = f(S_T)$, where $f$ is given in terms of a finite complex measure $\Pi$ [cf. (2)]. The choice of the set of *admissible trading strategies* is a delicate point in continuous time. Following [50], Section 1, we choose the set

$$\Theta := \left\{ \theta \in L(S) : \int_0^\cdot \theta_t \, dS_t \in \mathcal{H}^2 \right\},$$

which is well suited for quadratic optimization problems. Here, the space $\mathcal{H}^2$ of semimartingales is defined as follows:

DEFINITION 3.1. For any real-valued special semimartingale $Y$ with canonical decomposition $Y = Y_0 + N + B$, we define

$$\|Y\|_{\mathcal{H}^2} := \|Y_0\|_2 + \left\| \sqrt{[N,N]_T} \right\|_2 + \|\operatorname{var}(B)_T\|_2,$$

where $\operatorname{var}(B)$ denotes the variation process of $B$ and $\|\cdot\|_2$ denotes the $L^2$-norm. By $\mathcal{H}^2$ we denote the set of all real-valued special semimartingales $Y$ with $\|Y\|_{\mathcal{H}^2} < \infty$.

In our setup, this set can be expressed more easily as follows:

LEMMA 3.1. *We have*
$$\Theta = \left\{ \theta \text{ predictable process: } E\left(\int_0^T \theta_t^2 S_{t-}^2 \, dt\right) < \infty \right\}.$$

PROOF. From Lemma 3.2 below we conclude that $A_t = \kappa(1) \int_0^t S_{u-} \, du$ and

(13) $$\langle M, M \rangle_t = (\kappa(2) - 2\kappa(1)) \int_0^t S_{u-}^2 \, du$$

for the canonical decomposition $S = S_0 + M + A$ of the special semimartingale $S$. Hence we have

(14) $$A_t = \int_0^t \lambda_u \, d\langle M, M \rangle_u$$

with $\lambda_u := \lambda/S_{u-}$ and $\lambda := \frac{\kappa(1)}{\kappa(2) - 2\kappa(1)}$. Therefore, the mean-variance trade-off process

$$K_t = \int_0^t \lambda_u^2 \, d\langle M, M \rangle_u = \frac{\kappa(1)^2}{\kappa(2) - 2\kappa(1)} t$$

in the sense of [50], Section 1, is deterministic and bounded. According to [50], Lemma 2, we have that $\theta \in \Theta$ holds if and only if $\theta$ is predictable and $E(\int_0^T \theta_t^2 \, d\langle M, M \rangle_t) < \infty$. Since

$$\int_0^T \theta_t^2 \, d\langle M, M \rangle_t = (\kappa(2) - 2\kappa(1)) \int_0^T \theta_t^2 S_{t-}^2 \, dt,$$

the assertion follows. □

If we denote by $G_t(\theta) := \int_0^t \theta_u \, dS_u$ the *cumulative gains process* of $\theta \in \Theta$, then the *variance-optimal initial capital* and *variance-optimal hedging strategy* can be defined as in the previous section (with $T$ instead of $N$).

The following characterizations of the variance-optimal hedge and its expected squared error correspond to Theorems 2.1 and 2.2. They are proved at the end of this subsection.

THEOREM 3.1. *The variance-optimal initial capital $V_0$ and the variance-optimal hedging strategy $\phi$ are given by*
$$V_0 = H_0$$
*and the expression*

(15) $$\phi_t = \xi_t + \frac{\lambda}{S_{t-}}(H_{t-} - V_0 - G_{t-}(\phi)),$$

*where the processes $H_t$ and $\xi_t$ and the constant $\lambda$ are defined by*

$$\gamma(z) := \frac{\kappa(z+1) - \kappa(z) - \kappa(1)}{\kappa(2) - 2\kappa(1)},$$

$$\eta(z) := \kappa(z) - \kappa(1)\gamma(z),$$

(16) $$\lambda := \frac{\kappa(1)}{\kappa(2) - 2\kappa(1)},$$

$$H_t := \int S_t^z e^{\eta(z)(T-t)} \Pi(dz),$$

$$\xi_t := \int S_{t-}^{z-1} \gamma(z) e^{\eta(z)(T-t)} \Pi(dz).$$

*The optimal initial capital is unique. The optimal hedging strategy $\phi_t(\omega)$ is unique up to some $(P(d\omega) \otimes dt)$ null set.*

Remark 2.2 on risk-minimizing hedging for fixed initial endowment $c$ applies in continuous time as well.

THEOREM 3.2. *The variance of the hedging error $E((V_0 + G_T(\phi) - H)^2)$ in Theorem 3.1 equals*

$$J_0 := \int\int J_0(y,z) \Pi(dy) \Pi(dz),$$

*where*

$$\alpha(y,z) := \eta(y) + \eta(z) - \frac{\kappa(1)^2}{\kappa(2) - 2\kappa(1)},$$

$$\beta(y,z) := \kappa(y+z) - \kappa(y) - \kappa(z)$$
$$- \frac{(\kappa(y+1) - \kappa(y) - \kappa(1))(\kappa(z+1) - \kappa(z) - \kappa(1))}{\kappa(2) - 2\kappa(1)},$$

$$J_0(y,z) := \begin{cases} S_0^{y+z} \beta(y,z) \dfrac{e^{\alpha(y,z)T} - e^{\kappa(y+z)T}}{\alpha(y,z) - \kappa(y+z)}, & \text{if } \alpha(y,z) \neq \kappa(y+z), \\ S_0^{y+z} \beta(y,z) T e^{\kappa(y+z)T}, & \text{if } \alpha(y,z) = \kappa(y+z). \end{cases}$$

REMARK 3.1. If $(\mu, \sigma^2, \nu)$ denotes the Lévy–Khinchine triplet of $X$ (relative to the truncation function $x \mapsto x\mathbb{1}_{\{|x| \leq 1\}}$), then we have

$$\kappa(z) = \mu z + \frac{\sigma^2}{2} z^2 + \int (e^{zx} - 1 - zx\mathbb{1}_{\{|x| \leq 1\}}) \nu(dx)$$

for $z \in D$ (cf. [46], Theorem 25.17). In particular, we have $\kappa(z) = \mu z + (\sigma^2/2) z^2$ for Brownian motion. Note that

$$\Phi\left(\frac{x-\mu}{\sigma}\right) = \frac{1}{2\pi i} \int_{R-i\infty}^{R+i\infty} \frac{e^{(x-\mu)z + (\sigma^2/2)z^2}}{z} dz$$

for any $R > 0$, where $\Phi$ denotes the cumulative distribution function of $N(0,1)$. Using the same decomposition and substitution as in Remark 4.1, one easily shows that $V_0$ and $\phi$ in Theorem 3.1 coincide with the Black–Scholes price and the replicating strategy in the case of a European call $H$ and Brownian motion $X$. This does not come as a surprise because perfect hedging is clearly variance-optimal.

The hedging strategy $\phi$ in Theorem 3.1 is given in *feedback* form, that is, it is only known in terms of its own gains from trade up to time $t$. From a practical point of view, these gains are obviously known to the trader. However, they cannot be computed recursively as in the discrete-time case. Therefore, one may prefer an explicit expression for $G_t(\phi)$ from a mathematical point of view. It is provided by the following theorem.

THEOREM 3.3. *Suppose that $P(\Delta X_t = \log(1+1/\lambda)$ for some $t \in (0,T]) = 0$. Then the gains process of the variance-optimal hedging strategy $\phi$ in Theorem 3.1 is of the form*

$$G_t(\phi) = \mathcal{E}(-\lambda \widetilde{X})_t \left( \int_0^t \frac{\xi_u S_{u-} - \lambda(H_{u-} - V_0)}{\mathcal{E}(-\lambda \widetilde{X})_{u-}} \, dY_u \right),$$

*where the processes $\widetilde{X}$ and $Y$ are defined as*

(17)
$$\widetilde{X}_t := \mathcal{L}(S)_t := \int_0^t \frac{1}{S_{u-}} \, dS_u,$$
$$Y_t := \widetilde{X}_t + \int_0^t \frac{\lambda}{1 - \lambda \Delta \widetilde{X}_u} \, d[\widetilde{X}, \widetilde{X}]_u.$$

REMARK 3.2. The condition on $X$ is equivalent to assuming that the Lévy measure of $X$ puts no mass on $\log(1+1/\lambda)$. This holds for any model of practical importance. In the general case it is still possible to give an explicit, albeit more involved, representation of the gains process (cf. [28], (6.8)).

Moreover, observe that $\widetilde{X}$ and $Y$ are both Lévy processes (cf. [32], Lemma 2.7 and straightforward calculations). Recall that the stochastic exponential $\mathcal{E}(U)$ of a real-valued Lévy process or any other semimartingale $U$ can be written explicitly as

$$\mathcal{E}(U)_t = \exp(U_t - \tfrac{1}{2}[U,U]_t) \prod_{u \leq t} (1 + \Delta U_u) \exp(-\Delta U_u + \tfrac{1}{2}(\Delta U_u)^2)$$

(cf. [41], Theorem II.36).

The remainder of this section is devoted to the proof of Theorems 3.1–3.3. The approach parallels the one in the previous section. As before, we determine the Föllmer–Schweizer decomposition of the claim and apply results that relate this decomposition to the variance-optimal hedge.

LEMMA 3.2. *Let $z \in \mathbb{C}$ with $S_T^z \in L^2(P)$. Then $S^z$ is a special semimartingale whose canonical decomposition $S_t^z = S_0^z + M(z)_t + A(z)_t$ satisfies*

$$A(z)_t = \kappa(z) \int_0^t S_{u-}^z \, du \tag{18}$$

*and*

$$\langle M(z), M \rangle_t = (\kappa(z+1) - \kappa(z) - \kappa(1)) \int_0^t S_{u-}^{z+1} \, du, \tag{19}$$

*where $M = M(1)$ corresponds to $z = 1$ as in the proof of Lemma 3.1.*

PROOF. Note that almost by definition of the cumulant generating function, $N(z)_t := e^{-\kappa(z)t} S_t^z$ is a martingale. Integration by parts yields $S_t^z = e^{\kappa(z)t} N(z)_t = S_0^z + M(z)_t + A(z)_t$ with $M(z)_t = \int_0^t e^{\kappa(z)s} \, dN(z)_u$ and $A(z)$ as claimed. Moreover, we have

$$\begin{aligned}
[M(z), M]_t &= [S^z, S]_t \\
&= S_t^{z+1} - S_0^{z+1} - \int_0^t S_{u-}^z \, dS_u - \int_0^t S_{u-} \, dS_u^z \\
&= M(z+1)_t - \int_0^t S_{u-}^z \, dM_u - \int_0^t S_{u-} \, dM(z)_u \\
&\quad + (\kappa(z+1) - \kappa(z) - \kappa(1)) \int_0^t S_{u-}^{z+1} \, du.
\end{aligned}$$

Note that the first three terms on the right-hand side are local martingales. Since $\langle M(z), M \rangle$ is the predictable part of finite variation of the special semimartingale $[M(z), M]$, equation (19) follows. $\square$

DEFINITION 3.2. Denote by $S = S_0 + M + A$ the canonical special semimartingale decomposition of $S$ into a local martingale $M$ and a predictable process of finite variation $A$. The sum $H = H_0 + \int_0^T \xi_t \, dS_t + L_T$ is called a *Föllmer–Schweizer decomposition* of $H \in L^2(P)$ if $H_0$ is $\mathcal{F}_0$-measurable, $\xi \in \Theta$ and $L$ is a square-integrable martingale with $L_0 = 0$ that is orthogonal to $M$ (in the sense that $LM$ is a local martingale). We will use this notion as well if $H, H_0, \xi$ and $L$ are complex-valued, in which case we require $\Re(\xi) \in \Theta$ and $\Im(\xi) \in \Theta$.

The existence of a Föllmer–Schweizer decomposition was established in [50], Theorem 15, in our case of bounded mean-variance trade-off. It can be expressed in terms of a backward stochastic differential equation. Since the latter may be hard to solve, we do not use this result. Instead, we prove directly that the continuous-time limit of the expressions in Section 2 leads to a Föllmer–Schweizer decomposition.

LEMMA 3.3. *Let $z \in \mathbb{C}$ with $S_T^z \in L^2(P)$. Then $H(z) = S_T^z$ admits a Föllmer–Schweizer decomposition $H(z) = H(z)_0 + \int_0^T \xi(z)_t \, dS_t + L(z)_T$, where*

$$H(z)_t := e^{\eta(z)(T-t)} S_t^z,$$
(20)
$$\xi(z)_t := \gamma(z) e^{\eta(z)(T-t)} S_{t-}^{z-1},$$
$$L(z)_t := H(z)_t - H(z)_0 - \int_0^t \xi(z)_u \, dS_u$$

*and $\gamma(z), \eta(z)$ are defined in Theorem 3.1. Moreover, $M$ is a square-integrable martingale and hence $L(z)M$ is a martingale.*

PROOF. Partial integration and (18) yield

$$H(z)_t = H(z)_0 + \int_0^t e^{\eta(z)(T-u)} \, dM(z)_u + (\kappa(z) - \eta(z)) \int_0^t e^{\eta(z)(T-u)} S_{u-}^z \, du$$

and

$$\int_0^t \xi(z)_u \, dS_u = \int_0^t \xi(z)_u \, dM_u + \kappa(1)\gamma(z) \int_0^t e^{\eta(z)(T-u)} S_{u-}^z \, du.$$

Since $\kappa(z) - \eta(z) - \kappa(1)\gamma(z) = 0$, the predictable part of finite variation in the special semimartingale decomposition of $L(z)$ vanishes and we have

(21) $$L(z)_t = \int_0^t e^{\eta(z)(T-u)} \, dM(z)_u - \int_0^t \xi(z)_u \, dM_u,$$

which implies that $L(z)$ is a local martingale.

From (19) for $z$ and 1 instead of $z$ it follows that

$$\langle L(z), M \rangle_t = \int_0^t e^{\eta(z)(T-u)} \, d\langle M(z), M \rangle_u - \int_0^t \xi(z)_u \, d\langle M, M \rangle_u$$
$$= (\kappa(z+1) - \kappa(z) - \kappa(1) - \gamma(z)(\kappa(2) - 2\kappa(1)))$$
$$\quad \times \int_0^t e^{\eta(z)(T-u)} S_{u-}^{z+1} \, du$$
$$= 0.$$

Consequently, $L(z)M$ is a local martingale as well.

Similar calculations yield

$$\langle L(z), \overline{L(z)} \rangle_t = \langle L(z), L(\overline{z}) \rangle_t$$

(22)
$$= \left( \kappa(2\Re(z)) - 2\Re(\kappa(z)) - \frac{|\kappa(z+1) - \kappa(z) - \kappa(1)|^2}{\kappa(2) - 2\kappa(1)} \right)$$
$$\times \int_0^t e^{2\Re(\eta(z))(T-u)} S_{u-}^{2\Re(z)} \, du$$

and

(23)
$$\int_0^T |\xi(z)_t|^2 S_{t-}^2 \, dt = \left| \frac{\kappa(z+1) - \kappa(z) - \kappa(1)}{\kappa(2) - 2\kappa(1)} \right|^2$$
$$\times \int_0^T e^{2\Re(\eta(z))(T-t)} S_{t-}^{2\Re(z)} \, dt.$$

Since

(24)
$$E(S_{t-}^{2\Re(z)}) = E(S_t^{2\Re(z)}) = S_0^{2\Re(z)} e^{t\kappa(2\Re(z))} \leq S_0^{2\Re(z)} (1 \vee e^{T\kappa(2\Re(z))}) < \infty,$$

it follows that $E(\langle L(z), \overline{L(z)} \rangle_T) < \infty$. Therefore, $L$ is a square-integrable martingale.

Similarly, (23) and (24) yield that $\Re(\xi(z)) \in \Theta$ and $\Im(\xi(z)) \in \Theta$. Equations (19) and (24) for 1 instead of $z$ imply that $M$ is a square-integrable martingale. $\square$

LEMMA 3.4. *There exist constants $c_1, \ldots, c_5 \geq 0$ such that*

(25) $\quad \Re(\eta(z)) \leq c_1,$

$$0 \leq \kappa(2\Re(z)) - 2\Re(\kappa(z)) - \frac{|\kappa(z+1) - \kappa(z) - \kappa(1)|^2}{\kappa(2) - 2\kappa(1)}$$

(26)
$$\leq -c_2 \Re(\eta(z)) + c_3,$$
$$|\gamma(z)|^2 \leq -c_4 \Re(\eta(z)) + c_5$$

*for any $z \in \mathbb{C}$ with $R' \leq \Re(z) \leq R$, where $\gamma$ and $\eta$ are defined as in Theorem 3.1.*

PROOF. Since $\kappa$ is continuous, there is a constant $c_6 \geq 0$ such that

(27) $$|\kappa(2\Re(z))| \leq 2c_6$$

for any $z$ with $R' \leq \Re(z) \leq R$. Since $\langle L(z), \overline{L(z)} \rangle$ is increasing, (22) yields

$$\kappa(2\Re(z)) - 2\Re(\kappa(z)) - \frac{|\kappa(z+1) - \kappa(z) - \kappa(1)|^2}{\kappa(2) - 2\kappa(1)} \geq 0.$$

In particular
$$\Re(\kappa(z)) \leq \tfrac{1}{2}\kappa(2\Re(z)) \leq c_6$$

and

(28) $$\frac{|\kappa(z+1) - \kappa(z) - \kappa(1)|^2}{\kappa(2) - 2\kappa(1)} \leq 2c_6 - 2\Re(\kappa(z)),$$

which implies
$$\begin{aligned}|\kappa(1)\gamma(z)|^2 &\leq c_7 - c_8\Re(\kappa(z)) \\ &\leq c_9^2 + \tfrac{1}{4}(\Re(\kappa(z)))^2 \\ &\leq (|\tfrac{1}{2}\Re(\kappa(z))| + c_9)^2\end{aligned}$$

for some $c_7, c_8 \geq 0$ and $c_9 := \sqrt{c_7 + 4c_8^2}$. This yields

(29) $$\begin{aligned}\Re(\eta(z)) &= \Re(\kappa(z)) - \Re(\kappa(1)\gamma(z)) \\ &\leq \Re(\kappa(z)) + |\kappa(1)\gamma(z)| \\ &\leq c_{10} + \tfrac{1}{2}\Re(\kappa(z)) \\ &\leq c_9 + 2c_6 =: c_1\end{aligned}$$

with $c_{10} := c_9 + \tfrac{3}{2}c_6$. Inequality (28) also yields
$$|\gamma(z)|^2 \leq c_{11} - \frac{c_4}{2}\Re(\kappa(z))$$

for some $c_{11}, c_4 \geq 0$, which, together with (29), leads to
$$|\gamma(z)|^2 \leq c_{11} - c_4(\Re(\eta(z)) - c_{10}) = c_5 - c_4\Re(\eta(z))$$

with $c_5 := c_{11} + c_4 c_{10}$. Finally, the second inequality in (26) follows from (27), (29) and $\kappa(2) - 2\kappa(1) \geq 0$. □

PROPOSITION 3.1. *Any contingent claim $H = f(S_T)$ as in the beginning of this subsection admits a Föllmer–Schweizer decomposition $H = H_0 + \int_0^T \xi_t \, dS_t + L_T$. Using the notation of Lemma 3.3, it is given by*

(30) $$\begin{aligned}H_t &= \int H(z)_t \Pi(dz), \\ \xi_t &= \int \xi(z)_t \Pi(dz), \\ L_t &= \int L(z)_t \Pi(dz) = H_t - H_0 - \int_0^t \xi_u \, dS_u.\end{aligned}$$

*Moreover, the processes $H_t$, $\xi_t$ and $L_t$ are real-valued.*

PROOF. Let $z \in \mathbb{C}$ with $R' \leq \Re(z) \leq R$. Since $|H(z)_t|^2 = e^{2\Re(\eta(z))(T-t)} S_t^{2\Re(z)}$, we have that $E(|H(z)_t|^2)$ is bounded by some constant which depends only on $R$ and $R'$ [cf. (24) and (25)]. It follows that $H_t$ is a well-defined square-integrable random variable. Similarly, (22), (24) and Lemma 3.4 yield, after straightforward calculations, that

$$E(|L(z)_t|^2) = E(\langle L(z), \overline{L(z)} \rangle_t) \leq E(\langle L(z), \overline{L(z)} \rangle_T)$$

is bounded as well by such a constant. Therefore, $L_t$ is a well-defined square-integrable random variable as well. Finally, (23) and Lemma 3.4 yield that $E(|\xi(z)_t S_{t-}|^2)$ and also $E(\int_0^T |\xi(z)_u|^2 S_{u-}^2 \, du)$ are bounded by some constant which depends only on $t, R$ and $R'$. Therefore, $\xi$ is well defined, and $\Re(\xi) \in \Theta$ and $\Im(\xi) \in \Theta$ by Lemma 3.1. The same Fubini-type argument as in discrete time shows that $E((L_t - L_u)\mathbb{1}_B) = 0$ and $E((M_t L_t - M_u L_u)\mathbb{1}_B) = 0$ for $u \leq t$ and $B \in \mathcal{F}_u$ (cf. Proposition 2.1). Hence, $L$ is a square-integrable martingale which is orthogonal to $M$. To be precise, we interpret $L$ as the up to indistinguishability unique modification whose paths are almost surely cadlag (cf. [41], Corollary I.1). By Fubini's theorem for stochastic integrals (cf. [41], Theorem IV.46), we have

$$\int \int_0^t \xi(z)_u \, dS_u \Pi(dz) = \int_0^t \int \xi(z)_u \Pi(dz) \, dS_u = \int_0^t \xi_u \, dS_u.$$

Together with (30) and (20) it follows that $H_0, \xi$ and $L$ do indeed provide a Föllmer–Schweizer decomposition of $H$. As in the proof of Proposition 2.1, this time using [38], Theorem 3.4, instead of [48], the uniqueness of the real-valued Föllmer–Schweizer decomposition yields that the processes $H_t$, $\xi_t$ and $L_t$ are real-valued. □

PROOF OF THEOREM 3.1. According to the proof of Lemma 3.1, the mean-variance trade-off process of $S$ in the sense of [51], Section 1, equals

$$K_t = \frac{\kappa(1)^2}{\kappa(2) - 2\kappa(1)} t = \int_0^t \frac{\lambda}{S_{u-}} \, dA_u.$$

In view of Proposition 3.1, the optimality follows from Theorem 3 and Corollary 10 of [50].

As in the proof of Theorem 2.1, it follows that $V_0 = \widetilde{V}_0$ and $G_T(\phi) = G_T(\widetilde{\phi})$ if $\widetilde{V}_0$ and $\widetilde{\phi}$ denote another variance-optimal hedge. Observe that the local martingale $N_t := -\int_0^t \lambda_u \, dM_u$ satisfies $\langle N, N \rangle_T = \int_0^T \lambda_u^2 \, d\langle M, M \rangle_u = K_T$, where $\lambda_u$ is defined as in the proof of Lemma 3.1. From [5], Propositions 3.7 and 3.9(ii), and the remark after Definition 5.4, it follows that $G(\phi - \widetilde{\phi})$ is a $\mathcal{E}(N)$ martingale in the sense of that paper. By Proposition 3.12(i) in the same paper, it is determined by its terminal value $G_T(\phi - \widetilde{\phi}) = 0$, that is, $G_t(\phi - \widetilde{\phi}) = 0$ for any $t \in [0, T]$. Hence

$$\begin{aligned}
0 &= E([G(\phi - \widetilde{\phi}), G(\phi - \widetilde{\phi})]_T) \\
&= E\left(\int_0^T (\phi - \widetilde{\phi})_t^2 \, d[S, S]_t\right) \\
&= E\left(\int_0^T (\phi - \widetilde{\phi})_t^2 \, d[M, M]_t\right) \\
&= E\left(\int_0^T (\phi - \widetilde{\phi})_t^2 \, d\langle M, M\rangle_t\right) \\
&= (\kappa(2) - 2\kappa(1)) E\left(\int_0^T (\phi - \widetilde{\phi})_t^2 S_{t-}^2 \, dt\right).
\end{aligned}$$

This implies that $\phi_t(\omega) = \widetilde{\phi}_t(\omega)$ outside some $(P(d\omega) \otimes dt)$ null set. □

PROOF OF THEOREM 3.2. Similarly as in Lemma 3.2, it is shown that

$$\langle M(y), M(z)\rangle_t = (\kappa(y+z) - \kappa(y) - \kappa(z)) \int_0^t S_{u-}^{y+z} \, du.$$

From (21), $\langle L(y), M\rangle = 0$ and (19), it follows that

$$\begin{aligned}
\langle L(y), L(z)\rangle_t &= \int_0^t e^{(\eta(y)+\eta(z))(T-u)} \, d\langle M(y), M(z)\rangle_u \\
&\quad - \int_0^t \gamma(z) e^{(\eta(y)+\eta(z))(T-u)} S_{u-}^{z-1} \, d\langle M(y), M\rangle_u \\
&= \beta(y, z) \int_0^t e^{(\eta(y)+\eta(z))(T-u)} S_{u-}^{y+z} \, du.
\end{aligned}$$
(31)

Consequently,

$$\int_0^T e^{-(K_T - K_t)} \, d\langle L(y), L(z)\rangle_t = \beta(y, z) \int_0^T S_{t-}^{y+z} e^{\alpha(y,z)(T-t)} \, dt, \quad (32)$$

where $K$ denotes the mean-variance trade-off process as in the proof of Lemma 3.1. Since $E(S_{t-}^{y+z}) = S_0^{y+z} e^{\kappa(y+z)t}$, an application of Fubini's theorem yields

$$E\left(\int_0^T e^{-(K_T - K_t)} \, d\langle L(y), L(z)\rangle_t\right) = S_0^{y+z} \beta(y, z) \int_0^T e^{\kappa(y+z)t + \alpha(y,z)(T-t)} \, dt,$$

which equals $J_0(y, z)$.

Observe that

$$\Re\langle L(y), L(z)\rangle = \tfrac{1}{2}(\langle L(y) + L(\overline{z}), \overline{L(y) + L(\overline{z})}\rangle - \langle L(y), \overline{L(y)}\rangle - \langle L(z), \overline{L(z)}\rangle)$$

and

$$\langle L(y) + L(\overline{z}), \overline{L(y) + L(\overline{z})} \rangle$$
$$\leq \langle L(y) + L(\overline{z}), \overline{L(y) + L(\overline{z})} \rangle + \langle L(y) - L(\overline{z}), \overline{L(y) - L(\overline{z})} \rangle$$
$$= 2\langle L(y), \overline{L(y)} \rangle + 2\langle L(z), \overline{L(z)} \rangle.$$

In the proof of Proposition 3.1 we noted that $E(\langle L(z), \overline{L(z)} \rangle_T)$ and hence also the expected total variation of $\Re(\langle L(y), L(z) \rangle_t)$ is bounded by some constant which depends only on $R$ and $R'$. By replacing $L(\overline{z})$ with $iL(\overline{z})$, it follows analogously that the total variation of $\Im(\langle L(y), L(z) \rangle_t)$ is bounded by a similar constant. Therefore,

$$\int \int \langle L(y), L(z) \rangle_t \Pi(dy) \Pi(dz)$$

is a well-defined continuous, predictable, complex-valued process of finite variation. Since

$$L_t^2 = \int \int L(y)_t L(z)_t \Pi(dy) \Pi(dz),$$

an application of Fubini's theorem yields that

$$L_t^2 - \int \int \langle L(y), L(z) \rangle_t \Pi(dy) \Pi(dz)$$

is a martingale. This implies

$$\langle L, L \rangle_t = \int \int \langle L(y), L(z) \rangle_t \Pi(dy) \Pi(dz)$$

by definition of the predictable quadratic variation. Another application of Fubini's theorem yields

$$\int_0^T e^{-(K_T - K_t)} d\langle L, L \rangle_t = \int \int \int_0^T e^{-(K_T - K_t)} d\langle L(y), L(z) \rangle_t \Pi(dy) \Pi(dz)$$

and hence

$$E\left( \int_0^T e^{-(K_T - K_t)} d\langle L, L \rangle_t \right)$$
$$= \int \int E\left( \int_0^T e^{-(K_T - K_t)} d\langle L(y), L(z) \rangle_t \right) \Pi(dy) \Pi(dz)$$
$$= \int \int J_0(y, z) \Pi(dy) \Pi(dz).$$

By [50], Corollary 9, the left-hand side of the previous equation equals the variance of the hedging error. $\square$

Finally, we prove the explicit representation of the gains process.

PROOF OF THEOREM 3.3. By (15), $G(\phi)$ solves the stochastic differential equation

$$G_t(\phi) = \int_0^t \left(\xi_u + \frac{\lambda(H_{u-} - V_0)}{S_{u-}}\right) dS_u - \int_0^t \frac{\lambda}{S_{u-}} G_{u-}(\phi) \, dS_u$$

$$= \int_0^t (\xi_u S_{u-} + \lambda(H_{u-} - V_0)) \, d\widetilde{X}_u + \int_0^t G_{s-}(\phi) \, d(-\lambda \widetilde{X})_u.$$

By [28], (6.8), this equation has a unique solution, which is given by

$$G_t(\phi) = \mathcal{E}(-\lambda \widetilde{X})_t \bigg( \int_0^t \frac{\xi_u S_{u-} - \lambda(H_{u-} - V_0)}{\mathcal{E}(-\lambda \widetilde{X})_{u-}} \, d\widetilde{X}_u$$

$$+ \int_0^t \frac{\xi_u S_{u-} - \lambda(H_{u-} - V_0)}{\mathcal{E}(-\lambda \widetilde{X})_u} \, d[\widetilde{X}, \lambda \widetilde{X}]_u \bigg).$$

Since $\mathcal{E}(-\lambda \widetilde{X})_u = (1 - \lambda \Delta \widetilde{X}_u)\mathcal{E}(-\lambda \widetilde{X})_{u-}$, the assertion follows. $\square$

**4. Integral representation of payoff functions.** In the previous sections we used an integral representation of the payoff of the form (2). A precise characterization of the class of functions that allow for such a representation is delicate (cf. [6]) and of minor interest in this context. To derive hedging strategies and so forth, explicit formulas for $\Pi$ are required. They are provided here for a number of payoffs.

The basic example is, of course, the European call option $H = (S_N - K)^+$. Its integral representation (2) is provided by the following lemma.

LEMMA 4.1. *Let $K > 0$. For arbitrary $R > 1, s > 0$, we have*

$$(s - K)^+ = \frac{1}{2\pi i} \int_{R-i\infty}^{R+i\infty} s^z \frac{K^{1-z}}{z(z-1)} \, dz.$$

PROOF. For $\Re(z) > 1$, we have

$$\int_{-\infty}^{\infty} (e^x - K)^+ e^{-zx} \, dx = \frac{K^{1-z}}{z(z-1)}.$$

The assertion follows now from Theorem A.1. $\square$

REMARK 4.1. 1. Using $\frac{1}{z(z-1)} = \frac{1}{z-1} - \frac{1}{z}$ and substituting $z - 1$ for $z$ we can write the variance-optimal initial capital for the European call option as

$$V_0 = S_0 \Psi^{(1)}\left(\log\left(\frac{S_0}{K}\right)\right) - K \Psi^{(0)}\left(\log\left(\frac{S_0}{K}\right)\right)$$

with
$$\Psi^{(j)}(x) := \frac{1}{2\pi i} \int_{R-j-i\infty}^{R-j+i\infty} h(z+j)^N \frac{e^{zx}}{z} dz$$

in discrete time, respectively,
$$\Psi^{(j)}(x) := \frac{1}{2\pi i} \int_{R-j-i\infty}^{R-j+i\infty} \frac{e^{\eta(z+j)T+zx}}{z} dz$$

in continuous time. This resembles the pricing formulas for European calls in the Cox–Ross–Rubinstein, respectively, Black–Scholes, model. Note, however, that $\Psi^{(j)}(x)$ may not be a distribution function in general.

2. For the application of Lemma 4.1, we need slightly more than second moments of $X_1$ and hence $S_N$. This seems counterintuitive because the payoff grows only linearly in $S_N$. It is, in fact, possible to derive the optimal hedge in the case where only second moments exist. The idea is to consider the difference of the call and the stock [cf. (33)]. Since the stock itself corresponds to the unit mass $\Pi = \varepsilon_1$, one immediately obtains an integral representation (2) of the call in the strip $0 \leq \Re(z) \leq 1$. In this case, $\Pi$ is a complex measure concentrated on $\{1\} \cup (R + i\mathbb{R})$.

Let us consider the Laplace transform representations of some simpler payoff functions. They are mostly taken from [42] and they can be derived by straightforward calculations from Theorem A.1. Interestingly, the put option payoff is expressed by the same integral as the call, but with the vertical line of integration to the left of zero, that is,

$$(K-s)^+ = \frac{1}{2\pi i} \int_{R-i\infty}^{R+i\infty} s^z \frac{K^{1-z}}{z(z-1)} dz \qquad (R < 0).$$

A related example is the payoff

$$(33) \qquad (s-K)^+ - s = \frac{1}{2\pi i} \int_{R-i\infty}^{R+i\infty} s^z \frac{K^{1-z}}{z(z-1)} dz \qquad (0 < R < 1).$$

While this does not correspond to an option that arises in practice, it can be used to compute the variance-optimal hedge for calls and puts in a situation when the moment or cumulant function of the underlying asset exists in $0 \leq \Re(z) \leq 2$, but in no larger strip. This is actually the natural minimal integrability requirement in the present setup.

The *power call* (cf. [43]) can be represented by

$$((s-K)^+)^2 = \frac{1}{2\pi i} \int_{R-i\infty}^{R+i\infty} s^z \frac{2K^{2-z}}{z(z-1)(z-2)} dz \qquad (R > 2),$$

which generalizes to higher integer powers as

$$((s-K)^+)^n = \frac{1}{2\pi i} \int_{R-i\infty}^{R+i\infty} s^z \frac{n!K^{n-z}}{z(z-1)\cdots(z-n)} dz \qquad (R > n)$$

and even to arbitrary powers $\alpha > 1$ by

$$((s-K)^+)^\alpha = \frac{1}{2\pi i} \int_{R-i\infty}^{R+i\infty} s^z K^{\alpha-z} B(\alpha+1, z-\alpha)\, dz \qquad (R > \alpha),$$

where $B$ denotes the Euler beta function, which can be expressed by the more familiar Euler gamma function

$$B(\alpha, \beta) = \frac{\Gamma(\alpha)\Gamma(\beta)}{\Gamma(\alpha+\beta)}.$$

The *self-quanto call* can be written as

$$(s-K)^+ s = \frac{1}{2\pi i} \int_{R-i\infty}^{R+i\infty} s^z \frac{K^{1-z}}{(z-1)(z-2)}\, dz \qquad (R > 2).$$

The *digital option* with payoff function $f(s) = \mathbb{1}_{[K,\infty)}(s)$ coincides almost surely with the payoff function

(34) $$f(s) = \tfrac{1}{2}\mathbb{1}_{\{K\}}(s) + \mathbb{1}_{(K,\infty)}(s)$$

if the law of $S_N$, respectively, $S_T$, has no atoms. Using statement 2 in Theorem A.1, the latter can be expressed as

$$f(s) = \lim_{c\to\infty} \frac{1}{2\pi i} \int_{R-ic}^{R+ic} s^z \frac{K^{-z}}{z}\, dz \qquad (R > 0).$$

This suggests applying the results of the previous sections to the measure

(35) $$\Pi(dz) = \frac{1}{2\pi i} \frac{K^{-z}}{z}\, dz$$

in the case of the digital option. However, this measure is not of finite variation. Nevertheless, the main statements remain valid if we interpret the integrals as Cauchy principal value integrals.

LEMMA 4.2. *Theorems* 2.1, 2.2 *and* 3.1–3.3 *hold for the digital option* (34) *and the measure* (35) *if the integrals are interpreted in the principal value sense, that is,*

(36)
$$H_n := P\text{-}\lim_{c\to\infty} \int_{R-ic}^{R+ic} S_n^z h(z)^{N-n} \Pi(dz),$$

$$\xi_n := P\text{-}\lim_{c\to\infty} \int_{R-ic}^{R+ic} S_{n-1}^{z-1} g(z) h(z)^{N-n} \Pi(dz),$$

$$J_0 := \lim_{c\to\infty} \int_{R-ic}^{R+ic} \int_{R-ic}^{R+ic} J_0(y,z) \Pi(dy) \Pi(dz)$$

*and so forth, where* $P\text{-}\lim$ *refers to convergence in probability. In continuous time, the corresponding limit for* $\xi_t(\omega)$ *is to be interpreted in* $(P(d\omega) \otimes dt)$ *measure.*

Proof. We will show the assertion in the continuous-time setting. The discrete-time case follows similarly.

*Step* 1. For $c \in \mathbb{R}_+$ define $H^{(c)} := f^c(S_T)$ with

$$f^c(s) := \int_{R-ic}^{R+ic} s^z \Pi(dz).$$

Since $\overline{K^{-z}}/2\pi i z = -K^{-\overline{z}}/2\pi i \overline{z}$, it follows that $H^{(c)}$ is real-valued. For $s \in \mathbb{R}_+$, we have

$$f(s) - f^c(s) = \lim_{c' \to \infty} \frac{1}{2\pi i} \left( \int_{R+ic}^{R+ic'} \left(\frac{s}{K}\right)^z \frac{1}{z} dz + \int_{R-ic'}^{R-ic} \left(\frac{s}{K}\right)^z \frac{1}{z} dz \right)$$

$$= \lim_{c' \to \infty} \frac{1}{\pi} \int_c^{c'} \Re\left( \frac{(s/K)^{R+ix}}{R+ix} \right) dx.$$

The integrand equals

$$\left(\frac{s}{K}\right)^R \left( \frac{R\cos(x\log(s/K))}{R^2 + x^2} + \frac{R^2 \sin(x\log(s/K))}{(R^2 + x^2)x} - \frac{\sin(x\log(s/K))}{x} \right).$$

Since $\sup_{c \in \mathbb{R}_+} |\int_c^\infty \frac{\sin(x)}{x} dx| < \infty$ (cf. [1], Section 5.2), it follows that

$$\sup_{c \in \mathbb{R}_+} |f(s) - f^c(s)| \leq u s^R$$

for some $u \in \mathbb{R}_+$. Consequently, $(H^{(c)} - H)^2 \leq u^2 S_T^{2R} \in L^1$ for any $c \in \mathbb{R}_+$, which implies that $H^{(c)} \stackrel{c \to \infty}{\to} H$ in $L^2$ by dominated convergence.

*Step* 2. Denote by $H = \widetilde{H}_0 + \int_0^T \widetilde{\xi}_t dS_t + \widetilde{L}_T$ the Föllmer–Schweizer decomposition of $H$, which exists, for example, by [38], Theorem 3.4. Moreover, let $H_t^{(c)}$, $\xi_t^{(c)}$ and $L_t^{(c)}$ be defined as in Proposition 3.1 for the claim $H^{(c)}$. By Theorem 3.8 in [38], we have $H_0^{(c)} \to \widetilde{H}_0$,

$$(37) \qquad E\left( \int_0^T (\xi_t^{(c)} - \widetilde{\xi}_t)^2 d\langle M, M \rangle_t \right) \to 0$$

and $E((L_T^{(c)} - \widetilde{L}_T)^2) \to 0$ for $c \to \infty$. Since $L^{(c)}$ and $L$ are martingales, this implies $L_t^{(c)} \to \widetilde{L}_t$ in $L^2$ and hence in probability for any $t \in [0, T]$. Together with (14) we obtain

$$\int_0^t (\xi_u^{(c)} - \widetilde{\xi}_u) dM_u \to 0,$$

$$\int_0^t (\xi_u^{(c)} - \widetilde{\xi}_u) dA_u \to 0$$

and hence

$$\int_0^t \xi_u^{(c)} dS_u \to \int_0^t \widetilde{\xi}_u dS_u$$

in probability for any $t \in [0,T]$. Moreover, we have $\xi^{(c)} \to \widetilde{\xi}$ in measure relative to $P(d\omega) \otimes dt$ [cf. (37) and (13)]. Together we obtain that $\widetilde{H}_0, \widetilde{\xi}_t$ and $\widetilde{L}_t$ coincide with the expressions in Proposition 3.1 for $H$ if the integrals are interpreted in the principal value sense. Theorems 3.1 and 3.3 now follow precisely as in Section 3.

*Step* 3. Denote by $J_0^{(c)}$ and $\widetilde{J}_0$ the variance of the hedging error for $H^{(c)}$ and $H$, respectively. In a Hilbert space, the mapping $x \mapsto \|x - P(x)\|^2$ is continuous if $P$ denotes the projection on some closed subspace. Hence $J_0^{(c)} \to \widetilde{J}_0$ for $c \to \infty$. Since Theorem 3.2 is applicable to $H^{(c)}$, it follows that $\widetilde{J}_0$ coincides with $J_0$ in (36). □

The *log contract* of [39] does not seem to fit into this framework because the logarithm has no Laplace transform. Nevertheless we can express it as a difference of two payoffs, namely its positive and negative part. The former has a Laplace transform for $\Re(z) > 0$ and the latter for $\Re(z) < 0$, and we have

$$\log(s) = \frac{1}{2\pi i} \int_{R-i\infty}^{R+i\infty} s^z \frac{1}{z^2} \, dz - \frac{1}{2\pi i} \int_{R'-i\infty}^{R'+i\infty} s^z \frac{1}{z^2} \, dz$$

with $R' < 0$ and $R > 0$. In this case, $\Pi$ is a complex measure concentrated on $(R' + i\mathbb{R}) \cup (R + i\mathbb{R})$.

Finally, let us emphasize again that the whole approach is linear in the claim. Hence, we immediately obtain the variance-optimal hedge for any linear combination of the payoffs above as, for example, bull and bear spreads and collars.

**5. Examples with numerical illustrations.** In this section we illustrate how the approach is applied to concrete models that are considered in the literature. As an example, we provide numerical results for the normal inverse Gaussian model. The other setups lead to similar figures. Recall that all quantities are discounted, as in the theoretical developments above.

5.1. *Discrete-time hedging in the Black–Scholes model.* Suppose the underlying asset follows geometric Brownian motion with annual drift parameter $\mu$ and volatility $\sigma$. Then log returns per unit of time are normally distributed with mean $\mu - \sigma^2/2$ and variance $\sigma^2$.

Let us consider an option with maturity $T$ and a positive integer $N$. If trading is restricted to times $kT/N$ for $k = 0, 1, \ldots, N$, the market becomes incomplete. Theorem 2.1 applies with the moment generating function

$$m(z) = \exp\left(\left(\left(\mu - \frac{\sigma^2}{2}\right)z + \frac{\sigma^2 z^2}{2}\right)\frac{T}{N}\right).$$

If continuous trading is permitted, the Black–Scholes market is complete. Hence the hedging error is exactly zero. The variance-optimal capital and hedging strategy are given by the Black–Scholes price and delta hedging, respectively. It can be verified easily that this agrees in fact with the formulas in Theorems 3.1 and 3.2, where the relevant cumulant function is

$$\kappa(z) = \left(\mu - \frac{\sigma^2}{2}\right)z + \frac{\sigma^2 z^2}{2}.$$

Note that the drift parameter $\mu$ affects the discrete-time results but has no influence in continuous time.

5.2. *Merton's jump diffusion with normal jumps.* In the jump-diffusion model considered by [37], the logarithmic stock price is modeled as a Brownian motion with drift $\mu$ and volatility $\sigma$ plus occasional jumps from an independent compound Poisson process with intensity $\lambda$. A particularly simple and popular case is obtained when the jumps are normally distributed, say with mean $\nu$ and variance $\tau$. Then the moment function for discrete time is

$$m(z) = \exp\left(\left(\mu z + \frac{\sigma^2 z^2}{2} + \lambda(e^{\nu z + \tau^2 z^2/2} - 1)\right)\frac{T}{N}\right)$$

and the cumulant function for continuous time is

$$\kappa(z) = \mu z + \frac{\sigma^2 z^2}{2} + \lambda(e^{\nu z + \tau^2 z^2/2} - 1).$$

Note that Merton uses a slightly different parameterization.

It turns out that the variance-optimal capital for a European call option may be negative for some parameters, for example, for continuous-time hedging with $S_0 = 100$, $K = 110$, $T = 1$ and $\mu = 0.01$, $\sigma = 0.03$, $\lambda = 0.01$, $\nu = 0.2$, $\tau = 0.02$, we obtain $V_0 = -0.13$. This illustrates, as mentioned above, that, in general, the variance-optimal initial capital is not a price.

5.3. *Hyperbolic, normal inverse Gaussian and variance gamma models.* The hyperbolic, normal inverse Gaussian (NIG) and the variance gamma (VG) Lévy processes are subfamilies or limiting cases of the class of generalized hyperbolic models, which all fit in the general framework of this paper. We refer to [15] for further details.

5.3.1. *Hyperbolic model.* The moment generating function in the hyperbolic case is

$$(1) \qquad m(z) = \left(\frac{\sqrt{\alpha^2 - \beta^2}}{\sqrt{\alpha^2 - (\beta+z)^2}}\frac{K_1(\delta\sqrt{\alpha^2 - (\beta+z)^2})}{K_1(\delta\sqrt{\alpha^2 - \beta^2})}e^{\mu z}\right)^{T/N},$$

where $K_1$ denotes the modified Bessel function of the third kind with index 1. Some care has to be taken if $T/N$ is not an integer. The $T/N$th power in

(1) is, in fact, the distinguished $T/N$th power (see [46], Section 7). The cumulant function equals

$$\kappa(z) = \operatorname{Ln}\left(\frac{\sqrt{\alpha^2 - \beta^2}}{\sqrt{\alpha^2 - (\beta + z)^2}} \frac{K_1(\delta\sqrt{\alpha^2 - (\beta + z)^2})}{K_1(\delta\sqrt{\alpha^2 - \beta^2})} e^{\mu z}\right).$$

Here Ln denotes the distinguished logarithm; see [46], Section 7.

5.3.2. *Normal inverse Gaussian model.* The moment generating function of the normal inverse Gaussian model is given by

$$m(z) = \exp\left((\mu z + \delta(\sqrt{\alpha^2 - \beta^2} - \sqrt{\alpha^2 - (\beta + z)^2}))\frac{T}{N}\right).$$

Consequently, the cumulant function equals

$$\kappa(z) = \mu z + \delta(\sqrt{\alpha^2 - \beta^2} - \sqrt{\alpha^2 - (\beta + z)^2}).$$

5.3.3. *Variance gamma model.* If we subordinate a Brownian motion with drift $\beta$ by the gamma subordinator with parameters $\delta$ and $\alpha$, and add a linear drift with rate $\mu$, we obtain a variance gamma Lévy process. We refer to [35] for further details.

The corresponding moment function for discrete time is

$$m(z) = \left(e^{\mu z}\left(\frac{\alpha}{\alpha - \beta z - z^2/2}\right)^\delta\right)^{T/N}.$$

The cumulant function needed for continuous-time hedging is given by

$$\kappa(z) = \mu z + \delta \operatorname{Ln}\left(\frac{\alpha}{\alpha - \beta z - z^2/2}\right).$$

5.4. *Numerical illustration.* Figures 1–3 illustrate the results for a European call in the normal inverse Gaussian model, compared to Black–Scholes as a benchmark. The time unit is a year. The parameters of the normal inverse Gaussian distribution are $\alpha = 75.49$, $\beta = -4.089$, $\delta = 3.024$ and $\mu = -0.04$, which correspond to annualized values of the daily estimates from [45] for Deutsche Bank, assuming 252 trading days per year and discounting by an annual riskless interest rate of 4%.

The parameters for the benchmark Gaussian model are chosen such that both models lead to returns of the same mean and variance. We consider a European call option with strike price $K = 99$ and maturity $T = 0.25$, that is, three months from now.

Figure 1 shows the variance-optimal initial capital as a function of the stock price in the NIG model for both continuous and discrete time with $N = 12$ trading dates, that is, weekly rebalancing of the hedging portfolio.

The Black–Scholes price is plotted as well for comparison. One may observe that the three curves cannot be distinguished by eye, that is, they do not differ much in absolute terms. A similar picture is obtained in Figure 2 for the initial hedge ratio as a function of the initial stock price. The Black–Scholes delta provides a good proxy for the optimal hedge in the NIG model for both continuous and weekly rebalancing. As a result one may say that the Black–Scholes approach produces a reasonable hedge for the European call even if real data follow this rather different jump-type model.

The similarity ceases to hold when it comes to the hedging error, which vanishes in a true Black–Scholes world. Figure 3 shows the variance of the hedging error as a function of the number of trades, going from $N = 1$ (static hedging) to $N = 63$ (daily rebalancing). The horizontal line in Figure 3 indicates the variance of the hedging error for continuous rebalancing in the NIG model. The two decreasing curves refer to the discrete hedging error in the NIG and the Gaussian case, respectively. In the latter case it converges to 0, which is the error in the limiting Black–Scholes model. As far as the size is concerned, the variance of the error for $N = 12$ trading dates (weekly rebalancing) in the NIG setup ($1.04 \approx 1.02^2$) equals roughly the sum of the error in the corresponding Gaussian model (0.83, due to discrete rather than perfect hedging) and the inherent error in the continuous-time NIG model (0.257, due to incompleteness from jumps). The standard deviation 1.02 of the hedging error in the discrete NIG case may be compared to the Black–Scholes price 4.50 of the option.

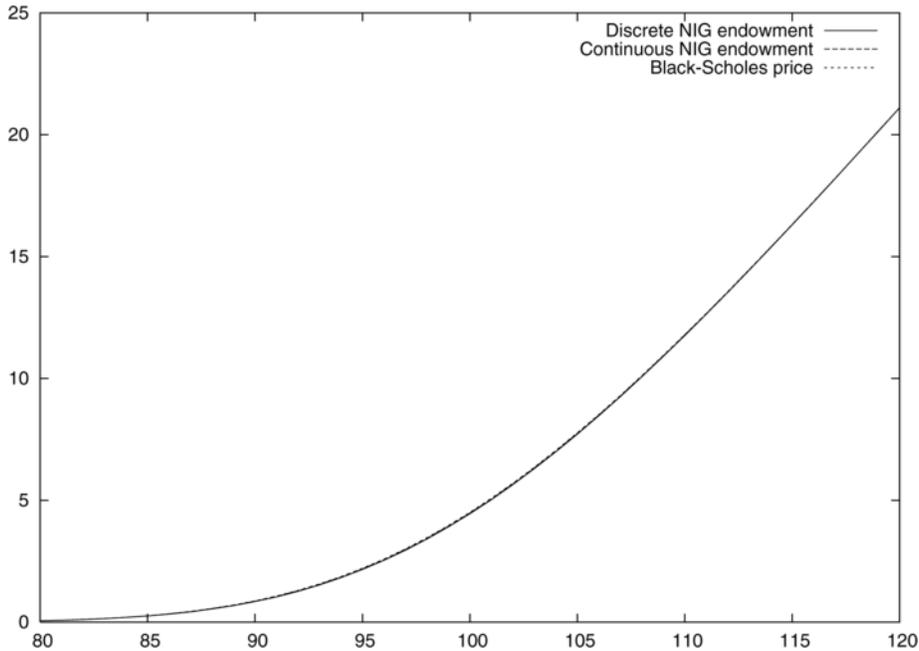

Fig. 1. *Variance-optimal initial capital for normal inverse Gaussian returns.*

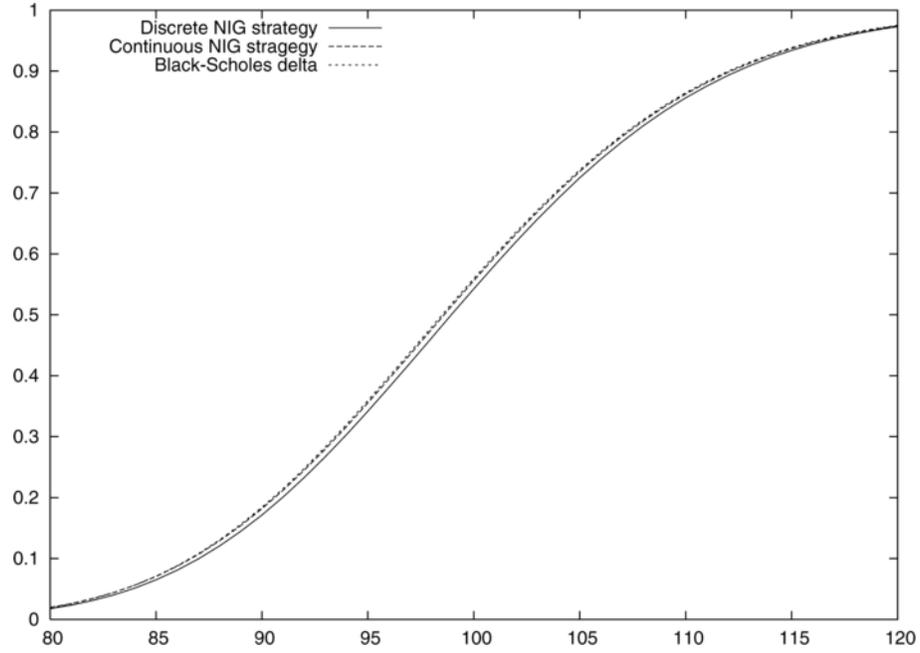

Fig. 2. *Variance-optimal initial hedge for normal inverse Gaussian returns.*

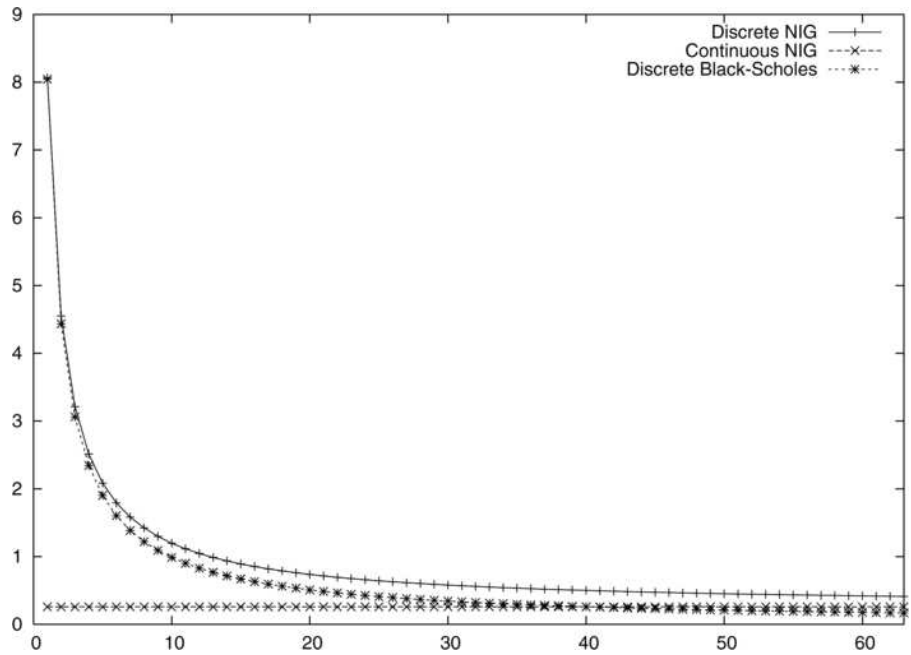

Fig. 3. *Variance of the hedging error for normal inverse Gaussian returns.*

## APPENDIX: BILATERAL LAPLACE TRANSFORMS

DEFINITION A.1. Let $f:\mathbb{R} \to \mathbb{C}$ be a measurable function. The (*bilateral*) *Laplace transform* $\widetilde{f}$ is given by

$$\text{(A.1)} \qquad \widetilde{f}(z) = \int_{-\infty}^{+\infty} f(x) e^{-zx}\, dx$$

for any $z \in \mathbb{C}$ such that the integral exists.

In the literature, the Laplace transform $\widetilde{f}$ is also denoted by $\mathcal{L}[f(x); z]$ or by $\mathcal{L}_{\mathrm{II}}[f(x); z]$ when it is necessary to distinguish the bilateral from the usual (unilateral) Laplace transform. The latter is defined by the same integral, but starting from 0 instead of $-\infty$.

We say that the Laplace transform $\widetilde{f}(z)$ exists if the Laplace transform integral (A.1) converges absolutely, or, in other words, if it exists as a proper Lebesgue integral as opposed to an improper integral. The following lemma shows that the domain of a Laplace transform is always a vertical strip in the complex plane. It may be empty, degenerate to a vertical line, a closed or open left or right half-plane, or all of $\mathbb{C}$.

LEMMA A.1. *Suppose that $\widetilde{f}(a)$ and $\widetilde{f}(b)$ exist for real numbers $a \leq b$. Then $\widetilde{f}(z)$ exists for any $z \in \mathbb{C}$ with $a \leq \Re(z) \leq b$.*

PROOF. The proof is obvious because $|f(x) e^{-zx}| = |f(x)| e^{-\Re(z) x} \leq |f(x) \times e^{-ax}| + |f(x) e^{-bx}|$. □

From

$$\text{(A.2)} \quad \widetilde{f}(u+iv) = \int_{-\infty}^{+\infty} f(x) e^{-(u+iv)x}\, dx = \int_{-\infty}^{+\infty} e^{ux} f(-x) e^{ixv}\, dx,$$

we see that $\mathcal{L}[f(x); u+iv] = \mathcal{F}[e^{ux} f(-x); v]$, where the last expression denotes the Fourier transform of the function $x \mapsto e^{ux} f(-x)$. Hence all properties of the bilateral Laplace transform can be reformulated in terms of the Fourier transform and vice versa.

There are many inversion formulas for the Laplace transform known in the literature. We will use the so-called Bromwich inversion integral, which can be justified by the following theorem.

THEOREM A.1. *Suppose that the Laplace transform $\widetilde{f}(R)$ exists for $R \in \mathbb{R}$.*

1. *If $v \mapsto \widetilde{f}(R + iv)$ is integrable, then $x \mapsto f(x)$ is continuous and*

$$f(x) = \frac{1}{2\pi i} \int_{R-i\infty}^{R+i\infty} \widetilde{f}(z) e^{zx} \, dz \qquad \text{for } x \in \mathbb{R}.$$

2. *If $f$ is of finite variation on any compact interval, then*

$$\lim_{\varepsilon \to 0} \frac{1}{2}(f(x+\varepsilon) + f(x-\varepsilon)) = \lim_{c \to \infty} \frac{1}{2\pi i} \int_{R-ic}^{R+ic} \widetilde{f}(z) e^{zx} \, dz \qquad \text{for } x \in \mathbb{R}.$$

PROOF. The first statement follows from [44], Theorem 9.11 and (A.2). For the second assertion, see [10], Satz 4.4.1. □

**Acknowledgments.** We would like to thank Arnd Pauwels for his comments and Thomas Jakob whose numerical computations lead to a correction of Figure 3. Thanks are also due two anonymous referees for valuable comments.

This paper is an extended and revised version of [27].

F. HUBALEK  
DEPARTMENT OF MATHEMATICAL SCIENCES  
UNIVERSITY OF AARHUS  
NY MUNKEGADE  
8000 ARHUS C  
DENMARK  
E-MAIL: fhubalek@imf.au.dk  

J. KALLSEN  
HVB-INSTITUTE FOR MATHEMATICAL FINANCE  
MUNICH UNIVERSITY OF TECHNOLOGY  
BOLTZMANNSTRASSE 3  
85747 GARCHING BEI MÜNCHEN  
GERMANY  
E-MAIL: kallsen@ma.tum.de

L. Krawczyk
Pl. Staromiejski 3/2
66-400 Gorzow Wlkp
Poland
E-mail: leszek_krawczyk@poczta.onet.pl